\documentclass[a4paper,11pt]{article}

\usepackage{amsmath,amsthm,amssymb}
\usepackage{amsfonts}
\usepackage{graphicx}
\usepackage{graphics}
\usepackage{url}
\usepackage{mathrsfs}
\usepackage{longtable}
\usepackage[dvips]{geometry}
\usepackage{algorithmic}
\usepackage{algorithm}
\usepackage[colorlinks=true,linkcolor=blue,citecolor=red]{hyperref}

\newtheorem{algo}{Algorithm}
\newenvironment{algthm}{\vskip2ex
\begin{algo}\rm\hrule\vskip1ex}
{\vskip1ex\nopagebreak\hrule\nopagebreak\end{algo}}

\oddsidemargin=.cm
\def\keywordname{{\bfseries Keywords}}
\def\keywords#1{\par\addvspace\medskipamount{\rightskip=0pt plus1cm
\def\and{\ifhmode\unskip\nobreak\fi\ $\cdot$
}\noindent\keywordname\enspace\ignorespaces#1\par}}%
\def\subclassname{{\bfseries Mathematics Subject Classification
(2020)}\enspace}
\def\subclass#1{\par\addvspace\medskipamount{\rightskip=0pt plus1cm
\def\and{\ifhmode\unskip\nobreak\fi\ $\cdot$
}\noindent\subclassname\ignorespaces#1\par}}
\pdfpagewidth=\paperwidth
\pdfpageheight=\paperheight

\title{Algorithms for pointwise and piecewise polynomial approximations to the trigonometric functions}
\author{LE PHUONG QUAN\thanks{
\parbox[t]{10cm}{\hbox{Department of Mathematics, College of Natural Sciences, Can Tho University}\vskip0.2ex
\hbox{3/2 Street, Can Tho City (900000), Vietnam}\vskip0.1ex
\hbox{Email:\,\url{lpquan@ctu.edu.vn}}}}}

\begin{document}
\maketitle
\begin{abstract}
In this paper, we propose a new and simple approach to the approximation algorithms that are modified and improved from our published results.
The computational and graphical examples are presented with the aid of Maple procedures.
\subclass{41A10 \and 41A80 \and 65D15 \and 65Y99}
\end{abstract}

\section{Introduction}\label{sect1}

We are mentioning the ancient problem: Finding polynomial approximations to the trigonometric functions. To propose a new solution to the problem, we will rely on some points in \cite{six}, especially on the method of keeping the integral part of $\dfrac{y}{(\pi/2)}$
unchanged when replacing $\pi/2$ with its approximate values for a given rational number $y\neq 0$. Accordingly, we will give a complete process with some modification for that method,
and present a different approach to our approximation algorithms.

Consider the function $f(y)=\sin y$. Given a real number $y\neq 0$ and an integer $n\ge 0$, we can write the Taylor formula for $f$ around $0$ as follows
\begin{equation}\label{eq1}
f(y)=\sum_{k=0}^n\frac{f^{(k)}(0)}{k!}y^k+\frac{f^{(n+1)}(\alpha)}{(n+1)!}y^{n+1},\quad \text{($\alpha$ is between $0$ and $y$)}.
\end{equation}
We know that $f^{(k)}(y)=\sin(y+k\pi/2)$, $k=0,1,\ldots,n$, so we have
$$f^{(k)}(0)=\sin(k\pi/2)=\left\{\begin{array}{@{}ccl}
0&&\text{as $k=2m\,\,(m=0,1,\ldots)$,}\\
(-1)^m&&\text{as $k=2m+1\,\,(m=0,1,\ldots)$.}
\end{array}\right.$$
Since $0\le k=2m+1\le n$, we infer the range for $m$: $0\le m\le\lfloor (n-1)/2\rfloor$, and we rewrite the sum in (\ref{eq1}) as
$$P_n(y)=\sum_{m=0}^{\lfloor (n-1)/2\rfloor}\frac{(-1)^m}{(2m+1)!}y^{2m+1}.$$
The polynomial $P_n(y)$ is called \textsl{Taylor polynomial\/} in the expansion (\ref{eq1}), and we can write
\begin{equation}\label{eq2}
\sin y=P_n(y)+\varepsilon_s,
\end{equation}
where
$$\varepsilon_s=\varepsilon_s(y,\alpha,n)=\frac{f^{(n+1)}(\alpha)}{(n+1)!}y^{n+1}.$$
If we choose the approximation $\sin y\approx P_n(y)$, then its absolute error satisfies
$$|\varepsilon_s|\le\frac{|y|^{n+1}}{(n+1)!}.$$
Similarly, for the function $g(y)=\cos y$, we can write
\begin{equation}\label{eq3}
\cos y=Q_n(y)+\varepsilon_c,
\end{equation}
where
$$Q_n(y)=\sum_{m=0}^{\lfloor n/2\rfloor}\frac{(-1)^m}{(2m)!}y^{2m},\quad\varepsilon_c=\varepsilon_c(y,\beta,n)=\frac{g^{(n+1)}(\beta)}{(n+1)!}y^{n+1}.$$
If we take the approximation $\cos y\approx Q_n(y)$, then its absolute error also satisfies
$$|\varepsilon_c|\le\frac{|y|^{n+1}}{(n+1)!}.$$

Such above approximations are meaningful if $|y|<1$ and $n$ may be chosen to be large enough to satisfy some given precision.
Note that the positions of $\alpha,\beta$ depend on arguments of $P_n,Q_n$. In the following, we will exploit these approximations effectively.

In general, for every real number $y$, there exists an integer $k$ such that
$$k\frac{\pi}{2}\le y<(k+1)\frac{\pi}{2}.$$
We will denote $\pi/2$ from now on by $p$, and we can write this inequality in the form of
$$k\le\frac{y}{p}<k+1,\quad k=\Big\lfloor\frac{y}{p}\Big\rfloor.$$
Now, we determine an integer $k_0$ by an if-then statement: if $y/p-\lfloor y/p\rfloor\le 1/2$, then $k_0=\lfloor y/p\rfloor$, else $k_0=\lfloor y/p\rfloor+1$.
In both cases, we always have $|y-k_0p|\le p/2<0.8$. This gives us a great advantage to approximate $\sin y$ in the following cases:
\begin{enumerate}
\item If $k_0$ is even, we can write
\begin{align*}
\sin y&=\sin(y-k_0p+k_0p)=\sin(y-k_0p)\cos(k_0p)+\sin(k_0p)\cos(y-k_0p)\\
&=\sin(y-k_0p)\cos(k_0p)=(-1)^{k_0/2}\sin(y-k_0p),
\end{align*}
then, from (\ref{eq2}), we can have the approximation
\begin{equation}\label{eq4}
\sin y\approx(-1)^{k_0/2}P_n(y-k_0p),
\end{equation}
and the absolute error $|\varepsilon|$ of this approximation satisfies
\begin{equation}\label{eq5}
|\varepsilon|\le\frac{|y-k_0p|^{n+1}}{(n+1)!}.
\end{equation}

\item If $k_0$ is odd, we can write
$$\sin y=\sin(k_0p)\cos(y-k_0p)=(-1)^{(k_0-1)/2}\cos(y-k_0p),$$
then, from (\ref{eq3}), we can have the approximation
\begin{equation}\label{eq6}
\sin y\approx(-1)^{(k_0-1)/2}Q_n(y-k_0p).
\end{equation}
The absolute error $|\varepsilon|$ of this approximation also satisfies (\ref{eq5}).
\end{enumerate}

In practice, we are only provided approximate values of $p$ and these values can be obtained from the lookup table-like storage of
approximate values of Pi (the number $\pi$). In some computer algebra systems, we can get the approximate value of Pi with $n$ digits by
their built-in commands such as \texttt{evalf[n](Pi)} with Maple, \texttt{vpa(pi,n)} with Matlab and \texttt{N[Pi,n]} with
Mathematica. As a convention here, we denote the approximate value of $p$ with $n$ digits by \texttt{valp[$n$]($p$)\/}.

In Section \ref{sect2}, before proposing a pointwise approximation algorithm, we will solve the two problems: how to choose approximate values $p'$ of $p$ such that $\lfloor y/p'\rfloor=\lfloor y/p\rfloor$
($y$ is a given rational number), and how to modify the approximations (\ref{eq4}) and (\ref{eq6}) when $p$ is replaced with $p'$ to guarantee the accuracy of $1/10^r$?
In Section \ref{sect3}, we propose a piecewise approximation algorithm based on the results obtained in Section \ref{sect2}. In Section \ref{sect4}, we present some other
approximation schemes together with estimates for their absolute error. In Section \ref{sect5}, we give some concluding remarks.

\section{Pointwise approximation algorithm}\label{sect2}

Firstly, we consider a rational number $y\neq 0$ and have the following relations
\begin{equation}\label{eq7}
0<\frac{y}{p}-\Big\lfloor\frac{y}{p}\Big\rfloor<1,\quad\frac{y}{p}-\Big\lfloor\frac{y}{p}\Big\rfloor\neq\frac{1}{2}
\end{equation}
and
\begin{equation}\label{eq8}
\max\Big\{\Big|\Big\lfloor\frac{y}{p}\Big\rfloor\Big|,\,\Big|\Big\lfloor\frac{y}{p}\Big\rfloor+1\Big|\Big\}<\frac{|y|}{1.5}+1.
\end{equation}

We will establish conditions on $p'$, an approximate value of $p$, such that
$$\Big\lfloor\frac{y}{p'}\Big\rfloor=\Big\lfloor\frac{y}{p}\Big\rfloor.$$
This relation is equivalent to
\begin{equation}\label{eq9}
\Big\lfloor\frac{y}{p'}\Big\rfloor<\frac{y}{p}<\Big\lfloor\frac{y}{p'}\Big\rfloor+1.
\end{equation}
We set $\sigma=|p-p'|$, then we can write
$$\Big|\frac{y}{p}-\frac{y}{p'}\Big|=\frac{\sigma|y|}{pp'}.$$
Now, (\ref{eq9}) is equivalent to
\begin{equation}\label{eq10}
\Big\lfloor\frac{y}{p'}\Big\rfloor-\frac{y}{p'}<\pm\frac{\sigma|y|}{pp'}<\Big\lfloor\frac{y}{p'}\Big\rfloor+1-\frac{y}{p'}.
\end{equation}
Since we may set $pp'>2.4$, we can have (\ref{eq10}) if $\sigma$ simultaneously satisfies the following conditions
\begin{align}
\sigma&<\frac{2.4}{|y|}\Big(\frac{y}{p'}-\Big\lfloor\frac{y}{p'}\Big\rfloor\Big)\label{eq11}\\
\sigma&<\frac{2.4}{|y|}\Big(\Big\lfloor\frac{y}{p'}\Big\rfloor+1-\frac{y}{p'}\Big).\label{eq12}
\end{align}
Suppose that we have both (\ref{eq11}) and (\ref{eq12}). We will check that the integer $k_0$ obtained in Section 1 is not changed when $p$ is replaced with $p'$.
Consider the following cases:
\begin{enumerate}
\item If $y/p'-\lfloor y/p'\rfloor<1/2$ and since
$$0<\frac{y}{p}-\Big\lfloor\frac{y}{p}\Big\rfloor=\frac{y}{p'}\pm\frac{\sigma|y|}{pp'}-\Big\lfloor\frac{y}{p'}\Big\rfloor,$$
we also have $y/p-\lfloor y/p\rfloor<1/2$ when $\sigma$ satisfies
\begin{equation}\label{eq13}
\sigma<\frac{2.4}{|y|}\Big(\frac{1}{2}+\Big\lfloor\frac{y}{p'}\Big\rfloor-\frac{y}{p'}\Big).
\end{equation}
Then, we choose $k_0=\lfloor y/p'\rfloor$ with this satisfaction.
\item If $y/p'-\lfloor y/p'\rfloor>1/2$, then we also have $y/p-\lfloor y/p\rfloor>1/2$ when $\sigma$ satisfies
\begin{equation}\label{eq14}
\sigma<\frac{2.4}{|y|}\Big(\frac{y}{p'}-\Big\lfloor\frac{y}{p'}\Big\rfloor-\frac{1}{2}\Big).
\end{equation}
Then, we choose $k_0=\lfloor y/p'\rfloor+1$ if we have (\ref{eq14}).
\item If $y/p'-\lfloor y/p'\rfloor=1/2$, then we also choose $k_0=\lfloor y/p'\rfloor$. Note that both (\ref{eq11}) and (\ref{eq12}) are obviously satisfied in this case when we
give an initial condition to $\sigma$ for all chosen values of $p'$:
\begin{equation}\label{eq15}
\sigma<\frac{1}{\Big(\dfrac{|y|}{1.5}+1\Big)10^{r+1}}.
\end{equation}
We refer to (\ref{eq8}) to motivate our choice of this condition.
\end{enumerate}

An important notice: (\ref{eq13}) implies (\ref{eq12}), and  (\ref{eq14}) implies (\ref{eq11}). We also emphasize that if $\sigma<1/10^m$, we can set $p'$ as the approximate value of
$p$ with the accuracy up to $m+2$ digits to prevent inexactitude from rounding-off rules. Therefore, when $m$ has been found such that
$\sigma<1/10^m$, we then set $p'=\texttt{valp[$m+2$]($p$)\/}$. We can find such an integer $m$ from the relation
\begin{equation}\label{eq16}
\frac{1}{10^m}\le\frac{1}{\Big(\dfrac{|y|}{1.5}+1\Big)10^{r+1}}<\frac{1}{10^{m-1}}.
\end{equation}
We will determine an appropriate value for $p'$ and then choose $k_0$ as mentioned above. For this task, we require $p'$ to satisfy both (\ref{eq11}) and (\ref{eq13}),
that is we may ask
\begin{equation}\label{eq17}
2.4\times 10^m\min\Big\{\frac{y}{p'}-\Big\lfloor\frac{y}{p'}\Big\rfloor,\,\frac{1}{2}+\Big\lfloor\frac{y}{p'}\Big\rfloor-\frac{y}{p'}\Big\}\ge|y|,
\end{equation}
or to satisfy both (\ref{eq12}) and (\ref{eq14}), that is we may ask
\begin{equation}\label{eq18}
2.4\times 10^m\min\Big\{\Big\lfloor\frac{y}{p'}\Big\rfloor-\frac{y}{p'}+1,\,\frac{y}{p'}-\Big\lfloor\frac{y}{p'}\Big\rfloor-\frac{1}{2}\Big\}\ge|y|.
\end{equation}
To prove the existence of such a $p'$, we need some results from the notion of sequence limit. This is really an important basis of our approximation algorithm.

Suppose that we have found $m=m_0$ satisfying (\ref{eq16}), then we choose $p_0=\texttt{valp[$m_0+2$]($p$)}$ and let $t_0=y/p_0-\lfloor y/p_0\rfloor$. We take an index $i$ that
starts with $0$. If $t_i=1/2$, then we choose $k_0=\lfloor y/p_i\rfloor$; otherwise, we consider the following process. From (\ref{eq17}) and (\ref{eq18}), we check
\begin{equation}\label{eq19}
2.4\times 10^{m_i}\min\{t_i,\,0.5-t_i\}\ge|y|,
\end{equation}
or
\begin{equation}\label{eq20}
2.4\times 10^{m_i}\min\{1-t_i,\,t_i-0.5\}\ge|y|.
\end{equation}
If both of (\ref{eq19}) and (\ref{eq20}) do not occur, we go to the next step by setting $m_{i+1}=m_i+1$, $p_{i+1}=\texttt{valp[$m_{i+1}+2$]($p$)}$ and
$t_{i+1}=y/p_{i+1}-\lfloor y/p_{i+1}\rfloor$. Then, we check again (\ref{eq19}) or (\ref{eq20}), but with $m_{i+1}$ and $t_{i+1}$ for this step.
Continuing this process, we go to the following cases:
\begin{enumerate}
\item One of either (\ref{eq19}) and (\ref{eq20}) first occurs at some step $k$, and assume it is (\ref{eq19}). Hence, we have $0<t_k=y/p_k-\lfloor y/p_k\rfloor<1/2$ (and also
$y/p-\lfloor y/p\rfloor<1/2$). Then, we choose $k_0=\lfloor y/p_k\rfloor$ ($=\lfloor y/p\rfloor$). Similarly, if (\ref{eq20}) occurs, we have
$1/2<t_k=y/p_k-\lfloor y/p_k\rfloor$ (and also $y/p-\lfloor y/p\rfloor>1/2$) and we choose $k_0=\lfloor y/p_k\rfloor+1=\lfloor y/p\rfloor+1$.
\item Both of (\ref{eq19}) and (\ref{eq20}) could not occur at any step. Then, we obtain the infinite sequences $m_i=m_0+i$, $p_i=\texttt{valp[$m_0+2+i$]($p$)}$ and
$t_i=y/p_i-\lfloor y/p_i\rfloor$, $i=0,1,\ldots$, such that
\begin{equation}\label{eq21}
2.4\times 10^{m_i}\min\{t_i,\,0.5-t_i\}<|y|,\quad\text{for all $i=0,1,\ldots$,}
\end{equation}
and
\begin{equation}\label{eq22}
2.4\times 10^{m_i}\min\{1-t_i,\,t_i-0.5\}<|y|,\quad\text{for all $i=0,1,\ldots$.}
\end{equation}
\end{enumerate}
Now, we prove the second case cannot be true. We have $p_i\to p$, then $y/p_i\to y/p$, as $i\to\infty$. Since $\lfloor y/p\rfloor<y/p<\lfloor y/p\rfloor+1$, there exists an
integer $i_0$ such that $\lfloor y/p\rfloor<y/p_i<\lfloor
y/p\rfloor+1$ for all $i\ge i_0$, hence $\lfloor
y/p_i\rfloor=\lfloor y/p\rfloor$ for all $i\ge i_0$. Therefore,
$$\lim_{i\to\infty}t_i=\lim_{i\to\infty}\Big(\frac{y}{p_i}-\Big\lfloor\frac{y}{p}\Big\rfloor\Big)=\frac{y}{p}-\Big\lfloor\frac{y}{p}\Big\rfloor.$$
If $y/p-\lfloor y/p\rfloor<0.5$, then
$$\min\{t_i,\,0.5-t_i\}\to\min\Big\{\frac{y}{p}-\Big\lfloor\frac{y}{p}\Big\rfloor,\,0.5-\frac{y}{p}+\Big\lfloor\frac{y}{p}\Big\rfloor\Big\}>0,$$
hence the left side of the inequality (\ref{eq21}) leads to $\infty$
as $i\to\infty$, and this is a contradiction. If $y/p-\lfloor y/p\rfloor>0.5$, by a similar argument, we also derive a
contradiction from (\ref{eq22}).

\begin{figure}[http]
\centering\includegraphics[width=9cm]{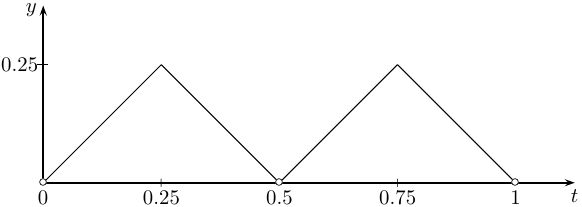}
\caption{The graph of $y=M(t)$ on
$(0,0.5)\cup(0.5,1)$.}\label{Fig1}
\end{figure}

There is one more question that should be considered explicitly.
What is the maximum value of $k$ that we can access the value \texttt{valp[$m_0+k+2$]($p$)}? Here, we provide an analysis
of this possibility. Let us take the sequence $\{t_k\}$, where $t_k=y/p_k-\lfloor y/p_k\rfloor$, $k=0,1,2,\ldots$, and $t_k\to y/p-\lfloor y/p\rfloor$ as $k\to\infty$.
Because $y\neq 0$ is rational, we must have (\ref{eq7}).
We set on the intervals $(0,0.5)$ and $(0.5,1)$ the functions $f(t)$ and $g(t)$, respectively, such that
$$f(t)=\min\{t,\,0.5-t\}=\left\{\begin{array}{@{}lcl}
t&&\text{as $0<t\le 0.25$,}\\
0.5-t&&\text{as $0.25<t<0.5$.}
\end{array}\right.$$
and
$${g(t)=\min\{1-t,t-0.5\}=\left\{\begin{array}{@{}lcl@{}}
t-0.5&&\text{as $0.5<t\le 0.75,$}\\
1-t&&\text{as $0.75<t<1.$}
\end{array}\right.}
$$
Then, we define
$$M(t)=\left\{\begin{array}{@{}lcl}
f(t)&&\text{as $0<t<0.5$,}\\
g(t)&&\text{as $0.5<t<1$.}
\end{array}\right.$$
and give its graph in Figure \ref{Fig1}.

Now, reaching (\ref{eq19}) or (\ref{eq20}) at the step $k$ can be expressed as
$$\frac{|y|}{2.4\times 10^{m_0}}\le 10^k M(t_k).$$
Because $|y|/10^{m_0}< 1/10^r$ (see (\ref{eq16}) where $m=m_0$), we can attain this when
\begin{equation}\label{eq23}
\frac{1}{2.4\times 10^r}\le 10^k M(t_k).
\end{equation}

In the case of $0<y/p-\lfloor y/p\rfloor<1/2$, we can conclude that there exist real numbers $\alpha$, $\beta$ and an integer $w_0$ such that
$$0<\alpha\le t_k\le\beta<\frac{1}{2}\quad
\text{for all $k\ge w_0$.}$$
Then, we easily obtain the relation
$$M(t_k)\ge\min\{M(\alpha),M(\beta)\}>0\quad\text{for all $k\ge w_0$.}$$
Hence, we can find the minimum value of $k$ satisfying (\ref{eq23}) or, equivalently, satisfying (\ref{eq19}) or (\ref{eq20}). We have the same argument if
$1/2<y/p-\lfloor y/p\rfloor<1$. In practice, since $|y|$ and $r$ are not too large, we can always access the value \texttt{valp[$m_0+k+2$]($p$)} with such a $k$.

Thus, for a given integer $r>0$, we can give a summary of the above process. At first, from an integer $m_0>0$ satisfying (\ref{eq16}), we choose $p_0=\texttt{valp[$m_0+2$]($p$)}$
and $t_0=y/p_0-\lfloor y/p_0\rfloor$. If $t_0=1/2$, we choose $k_0=\lfloor y/p_0\rfloor$. If $t_0\neq 1/2$, from the above diagram,
we can find an integer $k\ge 0$ and we choose $k_0=\lfloor y/p_k\rfloor$ or $k_0=\lfloor y/p_k\rfloor+1$, depending on $t_k=y/p_k-\lfloor y/p_k\rfloor$ and $m_k=m_0+k$ satisfy (\ref{eq19}) or (\ref{eq20})
($p_k=\texttt{valp[$m_0+k+2$]($p$)}$). Then, we take $n$ from (\ref{eq5}). Note that we always have $\lfloor y/p_k\rfloor=\lfloor y/p\rfloor$ ($k\ge 0$) and $|y-k_0p_k|\le p_k/2<0.8$.
This summary can be described more conveniently by the following pseudo-code algorithm.
\begin{algthm}\label{HQ-Algo1}
Determining $p'$, $k_0$ and $n$ from a given rational number $y\neq 0$.
\end{algthm}
\begin{algorithmic}[1]
\REQUIRE A rational number $y\neq 0$;
\ENSURE The values $p'$, $k_0$ and $n$;
\STATE $m:=r+1$;
\WHILE{$\dfrac{1}{\Big(\dfrac{|y|}{1.5}+1\Big)10^{r+1}}<\dfrac{1}{10^m}$}
\STATE $m:=m+1$; \COMMENT{Finding $m$ that satisfies (\ref{eq16})};
\ENDWHILE
\STATE $p':=\texttt{valp[$m+2$]($p$)}$, $t:=y/p'-\lfloor y/p'\rfloor$;
\IF{$t=0.5$}
\STATE $k_0:=\lfloor y/p'\rfloor$;
\ELSE
\WHILE{$m,t$ do not satisfy both (\ref{eq19}) and (\ref{eq20})}
\STATE $m:=m+1$, $p':=\texttt{valp[$m+2$]($p$)}$, $t:=y/p'-\lfloor y/p'\rfloor$;
\ENDWHILE
 \IF{$m,t$ satisfy (\ref{eq19}) (or $t<0.5$)}
 \STATE $k_0:=\lfloor y/p'\rfloor$;
 \ELSE
 \STATE $k_0:=\lfloor y/p'\rfloor+1$;
 \ENDIF
\ENDIF
\STATE $n:=0$;
\WHILE{$\dfrac{|y-k_0p'|^{n+1}}{(n+1)!}\ge\dfrac{1}{10^{r+1}}$}
\STATE $n:=n+1$; \COMMENT{Finding $n$ that satisfies $|\varepsilon|<\dfrac{1}{10^{r+1}}$ (from (\ref{eq5}))};
\ENDWHILE
\RETURN $p'$, $k_0$ and $n$;
\vskip1ex
\hrule
\vskip2ex
\end{algorithmic}

Finally, we will examine the accuracy of $1/10^r$ for the approximations (\ref{eq4}) and (\ref{eq6}) when we replace $p$ with its approximate values $p'$ obtained
from Algorithm \ref{HQ-Algo1}. We consider the following cases:
\begin{enumerate}
\item $y$ is a rational number and $y\neq 0$: From Algorithm \ref{HQ-Algo1}, we obtain the values $p'$, $k_0$ and $n$.
Assume that $k_0$ is even. Since $|y-k_0p'|<0.8$, we can approximate $\sin(y-k_0p')$ by the polynomial on the right-hand side of (\ref{eq4}) when $p$ is replaced by $p'$.
Therefore, the absolute error of this approximation is less than $1/10^{r+1}$. We denote this approximate value of $\sin(y-k_0p')$ by $A$. Then we can write
\begin{align*}
\sin y&=(-1)^{k_0/2}\sin(y-k_0p)\\
&=(-1)^{k_0/2}[\sin(y-k_0p)-\sin(y-k_0p')]+(-1)^{k_0/2}\sin(y-k_0p').
\end{align*}
Hence, we have the estimate
\begin{align*}
\big|\sin y-(-1)^{k_0/2}A\big|&\le|\sin(y-k_0p)-\sin(y-k_0p')|+|\sin(y-k_0p')-A|\\
&<|(y-k_0p)-(y-k_0p')|+\frac{1}{10^{r+1}}=|k_0|\sigma+\frac{1}{10^{r+1}}\\
&<\frac{|k_0|}{\Big(\dfrac{|y|}{1.5}+1\Big)10^{r+1}}+\frac{1}{10^{r+1}}<\frac{1}{10^{r+1}}+\frac{1}{10^{r+1}}<\frac{1}{10^r}.
\end{align*}
Now, assume that $k_0$ is odd. Then we can approximate $\sin(y-k_0p')$ by the polynomial on the right-hand side of (\ref{eq6}) when $p$ is replaced by $p'$.
If we denote this approximate value of $\sin(y-k_0p')$ by $B$ and use a similar argument as above, we can obtain the estimate
$$|\sin y-(-1)^{(k_0-1)/2}B|<\frac{1}{10^r}.$$
\item $y$ is an irrational number: Assume that we can find a rational number $y'\ne 0$ such that $|y-y'|<1/10^{r+1}$. According to the first case, we can approximate $\sin y'$ by a
number $C$ with the accuracy of $1/10^{r+1}$. Then, we have
$$|\sin y-C|=|\sin y-\sin y'+\sin y'-C|\le|y-y'|+|\sin y'-C|<\frac{1}{10^{r+1}}+\frac{1}{10^{r+1}}<\frac{1}{10^r}.$$
\end{enumerate}
There are many algorithms in the literature to find rational approximations to an irrational number with the desired accuracy, mostly using continued fractions.
The theoretical basis of this classic problem can be found in the two great books \cite{two} and \cite{five}. Recently, the survey article
\cite{one} applies a particular criterion to prove irrationality of a number and gives a nice approximation to it by a sequence of rational numbers.

We give an example on evaluating approximate values of $\sin(\pi/k)$ with the accuracy of $1/10^r$, where $k$ is an integer and $k\ge 2$. By choosing $p'=\texttt{valp[$r+3$]($p$)}$, we
can write
$$p=p'+\varepsilon,\quad 0<\varepsilon<\frac{1}{10^{r+1}}.$$
Hence, we have
$$\frac{\pi}{k}=\frac{2p}{k}=\frac{2}{k}p'+\frac{2}{k}\varepsilon,\quad\frac{2}{k}\varepsilon<\frac{1}{10^{r+1}}.$$
Therefore, $y=(2/k)\texttt{valp[$r+3$]($p$)}$ satisfies $|y-\pi/k|<1/10^{r+1}$. Based on Algorithm \ref{HQ-Algo1} and the two cases just mentioned above we can make a procedure to evaluate
the approximate values of $\sin(\pi/k)$ with the accuracy of $1/10^{10}$, $1/10^{20}$, $1/10^{30}$ and $1/10^{50}$. Some of these values are given in the following table:

\begin{longtable}[c]{|c|c|c|}\hline
$k$&$r$&$\approx\sin(\pi/k)$\\ \hline
\endfirsthead\hline
$k$&$r$&$\approx\sin(\pi/k)$\\ \hline
\endhead
\hline\multicolumn{3}{r}{\emph{continued on the next page}}
\endfoot
\hline
\endlastfoot
$3$&$10$&$0.8660254037$\\
$5$&$10$&$0.5877852524$\\
$7$&$20$&$0.43388373919609773406$\\
$9$&$20$&$0.34202014332673442399$\\
$17$&$30$&$0.183749517857635869327474437337$\\
$31$&$50$&$0.10116832198455048299912102931663774619952792876512$
\end{longtable}

\section{Piecewise approximation algorithm}\label{sect3}

In this section, we are interested in approximations to the sine function on an interval $[a,b]$. So we see Algorithm \ref{HQ-Algo1} as a tool that also supplies an
approximate polynomial from a given rational $y\neq 0$. More clearly, after applying Algorithm \ref{HQ-Algo1} to $y$, we obtain $p'$, $k_0$ and $n$; then,
we can approximate the value $\sin y$ by $A$ or $B$, depending on $k_0$ is even or odd. That is, we can write
$$\sin y\approx A=A(y)=(-1)^{k_0/2}P_n(y-k_0p')\quad\text{($k_0$ is even)}$$
or
$$\sin y\approx B=B(y)=(-1)^{(k_0-1)/2}Q_n(y-k_0p')\quad\text{($k_0$ is odd)},$$
where $|y-k_0p'|\le p'/2$ or $y\in[k_0p'-p'/2,k_0p'+p'/2]$ and $n$ is \textsl{now\/} the smallest integer such that
\begin{equation}\label{eq24}
\frac{(0.8)^{n+1}}{(n+1)!}<\frac{1}{10^{r+1}}.
\end{equation}
We choose such an $n$, which is the degree of $P_n$ or $Q_n$, to guarantee the accuracy of $1/10^r$ with the approximation $\sin x\approx A(x)$ or $\sin x\approx B(x)$ for
all $x\in[k_0p'-p'/2,k_0p'+p'/2]$. Because if we chose $n$ as in Algorithm \ref{HQ-Algo1}, there would be a $z\in[k_0p'-p'/2,k_0p'+p'/2]$ such  that $|z-k_0p'|>|y-k_0p'|$, and
the precision of $\sin z\approx A(z)$ or $\sin z\approx B(z)$ could not be $1/10^r$.
Therefore, for every \textsl{real number\/} $x\in[k_0p'-p'/2,k_0p'+p'/2]$ or $|x-k_0p'|\le p'/2$, by using arguments similar to the ones at the end of
Section \ref{sect2}, we can have the approximations $\sin x\approx A(x)$ or $\sin x\approx B(x)$ with the accuracy of $1/10^r$, depending on $k_0$ is even or odd. So, if we
denote both $A$ and $B$ by $F_y$, we can write
\begin{equation}\label{eq25}
\sin\approx\text{$F_y$ on $[(k_0-1/2)p',(k_0+1/2)p']$,}
\end{equation}
with the accuracy of $1/10^r$. We call $[(k_0-1/2)p',(k_0+1/2)p']$ \textsl{the approximate interval generated by $y$\/} or we call
$y$ \textsl{the generating point of the approximate interval\/} $[(k_0-1/2)p',(k_0+1/2)p']$. Note that for every real number $x$ such that $|x-k_0p'|\le 0.8$, we can also approximate
$\sin x$ by $A(x)$ or $B(x)$ with the accuracy of $1/10^r$. Thus, we may replace the interval in (\ref{eq25}) with $[k_0p'-0.8,k_0p'+0.8]$.

To emphasize the above process of getting the components $p'$, $k_0$ and $F_y$, corresponding to a given rational number $y\neq 0$, we denote it by
\texttt{IntvApprox} and present it as a map
$$\texttt{IntvApprox}\colon y\longmapsto[p',k_0,F_y].$$
To access the components of the list, we set
$$\texttt{IntvApprox($y$)[$1$]}=p',\quad \texttt{IntvApprox($y$)[$2$]}=k_0,\quad \texttt{IntvApprox($y$)[$3$]}=F_y.$$

Now, we go to the problem: Find a piecewise function $F$ that approximates the sine function
on an interval $[a,b]$ with the accuracy of $1/10^r$. Here, we can assume $a,b$ are rational, because we can always find rational numbers $\alpha,\beta$ such that $[a,b]\subset[\alpha,\beta]$.
Because the sine function is odd, we first determine the function $F$ on $[a,b]$ with $0\le a$. If $0\le a<b\le 0.8$, we can take $F=F_0$ on $[a,b]$,
where $F_0$ is given by:
$$F_0(x)= P_n(x)=\sum_{m=0}^{\lfloor(n-1 )/2\rfloor}\frac{(-1)^m}{(2m+1)!}x^{2m+1},$$
and $n$ is determined from (\ref{eq24}); hence, if $0\le a<0.8<b$, we set $F=F_0$ on $[a,0.8]$.
Therefore, it is sufficient to find $F$ on $[a,b]$ when $0.8\le a$. Before solving this problem in the following,
we set a global condition for $\sigma=|p-p'|$ as
\begin{equation}\label{eq26}
\sigma<\frac{1}{(2b+4)10^{r+1}}.
\end{equation}
This setting is the invocation of the estimates in our arguments.
We assert that $k<b+1.5$ for all integer $k$ obtained from generating points in our approximation algorithm. This assertion will be clearly explained later.
Also, we will check that these generating points satisfy the condition (\ref{eq15}).

\begin{figure}[htbp]
\centering
\includegraphics{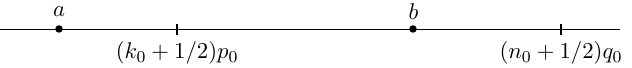}
\caption{The right ending points of the approximate intervals generated by $a$ and $b$.}\label{Fig2}
\end{figure}

First of all, we suppose that
$$\texttt{IntvApprox($a$)}=[p_0,k_0,A_0]\quad\text{and}\quad \texttt{IntvApprox($b$)}=[q_0,n_0,B_0].$$
These initial components are simply illustrated in Figure \ref{Fig2}, and we can see what should be done next, if we want to spread approximate intervals to the right.

It is easy to check that $k_0>0$ and if $k_0\ge n_0$, we must have $k_0-n_0\le 1$. We determine $F$ in the following cases:
\begin{enumerate}
\item If $k_0\ge n_0$, we claim $|b-k_0p_0|<0.8$. Indeed, we have
\begin{align*}
|b-k_0p_0|&=|b-n_0q_0+n_0q_0-k_0p_0|\le|b-n_0q_0|+|p_0-q_0|\\
&\le\frac{q_0}{2}+2\sigma<\frac{q_0}{2}+\frac{2}{(2b+4)10^{r+1}}<\frac{1.5708}{2}+\frac{1}{2\times10^{r+1}}\le 0.7854+0.005<0.8.
\end{align*}
Then, we choose $F=A_0$ on $[a,b]$.
\item If $k_0<n_0$, we take $(k_0+1)p_0$ as the next generating point; so, if
$$\texttt{IntvApprox($(k_0+1)p_0$)}=[p_1,k_1,A_1],$$
we must have $k_1=k_0+1\le n_0$, due to $|(k_0+1)p_0-k_1p_1|\le p_1/2$. Here, we choose this generating point to get the approximate interval
$[(k_1-1/2)p_1,(k_1+1/2)p_1]$ that can be seen as the consecutiveness of the interval $[(k_0-1/2)p_0,(k_0+1/2)p_0]$, since
$$(k_0+1/2)p_0=(k_1-1/2)p_0\approx(k_1-1/2)p_1$$
and $|(k_0+1/2)p_0-k_1p_1|<0.8$. Indeed, from (\ref{eq26}), we have
\begin{align*}
|(k_0+1/2)p_0-k_1p_1|&=|k_0(p_1-p_0)+p_1-p_0/2|\le 2k_0\sigma+\sigma+\frac{p_1}{2}\\
&<\frac{1}{10^{r+1}}+\frac{1}{5\times 10^{r+1}}+\frac{1.5708}{2}\le0.01+0.002+0.7854<0.8.
\end{align*}
Thus, we can set $F=A_1$ on $[(k_0+1/2)p_0,(k_0+1)p_0]$. If $k_1=n_0$, depending on where $b$ is from $(k_0+1)p_0$ (see Figure \ref{Fig3}), we set $F$ on $[a,b]$ as follows:
\begin{figure}[htbp]
\centering
\includegraphics{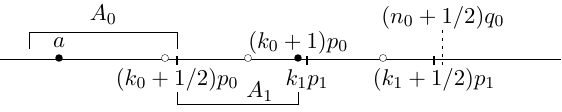}
\caption{The white nodes indicate the possible positions of $b$ in the case of $k_1=n_0$.}\label{Fig3}
\end{figure}
\begin{description}
\item[$b\le(k_0+1/2)p_0\colon$] $F=A_0$ on $[a,b]$.
\item[$(k_0+1/2)p_0<b\le(k_0+1)p_0\colon$]
$$F=\left\{\begin{array}{@{}lcl}
A_0&&\text{on $[a,(k_0+1/2)p_0)$,}\\
A_1&&\text{on $[(k_0+1/2)p_0,b]$.}
\end{array}
\right.$$
\item[$(k_0+1)p_0<b\colon$]
$$F=\left\{\begin{array}{@{}lcl}
A_0&&\text{on $[a,(k_0+1/2)p_0)$,}\\
A_1&&\text{on $[(k_0+1/2)p_0,(k_0+1)p_0)$,}\\
B_0&&\text{on $[(k_0+1)p_0,b]$.}
\end{array}
\right.$$
\end{description}
\item If $k_1<n_0$, we continue to choose $(k_1+1)p_1$ as the next generating point and get
$$\texttt{IntvApprox($(k_1+1)p_1$)}=[p_2,k_2,A_2].$$
We also imply $k_2=k_1+1\le n_0$ and $|(k_1+1/2)p_1-k_2p_2|<0.8$. Thus, we can set $F=A_1$ on $[(k_0+1/2)p_0,(k_1+1/2)p_1]$ and $F=A_2$ on
$[(k_1+1/2)p_1,(k_1+1)p_1]$. If $k_2=n_0$, depending on where $b$ is from $(k_1 + 1)p_1$ (see Figure \ref{Fig4}), we set $F$ on $[a,b]$ as follows:
\begin{figure}[htbp]
\centering
\includegraphics{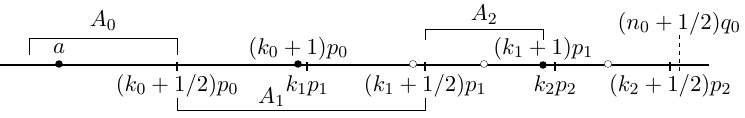}
\caption{The white nodes indicate the possible positions of $b$ in the case of $k_2=n_0$.}\label{Fig4}
\end{figure}
\begin{description}
\item[$b\le(k_1+1/2)p_1\colon$]
$$F=\left\{\begin{array}{@{}lcl}
A_0&&\text{on $[a,(k_0+1/2)p_0)$,}\\
A_1&&\text{on $[(k_0+1/2)p_0,b]$.}
\end{array}
\right.$$
\item[$(k_1+1/2)p_1<b\le(k_1+1)p_1\colon$]
$$F=\left\{\begin{array}{@{}lcl}
A_0&&\text{on $[a,(k_0+1/2)p_0)$,}\\
A_1&&\text{on $[(k_0+1/2)p_0,(k_1+1/2)p_1)$,}\\
A_2&&\text{on $[(k_1+1/2)p_1,b]$.}
\end{array}
\right.$$
\item[$(k_1+1)p_1<b\colon$]
\begin{equation}\label{eq27}
F=\left\{\begin{array}{@{}lcl}
A_0&&\text{on $[a,(k_0+1/2)p_0)$,}\\
A_1&&\text{on $[(k_0+1/2)p_0,(k_1+1/2)p_1)$,}\\
A_2&&\text{on $[(k_1+1/2)p_1,(k_1+1)p_1)$,}\\
B_0&&\text{on $[(k_1+1)p_1,b]$.}
\end{array}
\right.
\end{equation}
\end{description}
\end{enumerate}

By similar arguments as in the above cases, we finally obtain a finite sequence
$$\texttt{IntvApprox($(k_{i-1}+1)p_{i-1}$)}=[p_i,k_i,A_i],\quad i=1,2,\ldots,n,$$
where $k_i=k_{i-1}+1$, $i=1,2,\ldots,n-1$ and $k_n=n_0$. Then, we can check that
$$|(k_{i-1}+1/2)p_{i-1}-k_ip_i|<0.8,\quad i=1,2,\ldots,n,$$
so $[(k_i-1/2)p_i,(k_i+1/2)p_i]$ is consecutive to $[(k_{i-1}-1/2)p_{i-1},(k_{i-1}+1/2)p_{i-1}]$, $i=1,2,\ldots,n-1$.

Before finishing the process to determine the piecewise function $F$ on $[a,b]$ with $0.8\le a$, we give here
the explanation of the setting (\ref{eq26}). Since $|b-n_0q_0|\le q_0/2$, we get
$$n_0\le\frac{b}{q_0}+\frac{1}{2}<b+0.5,$$
so, we have that $k_i\le n_0+1<b+1.5$ for all $i=0,\ldots,n$. Then, it is easy to see that
$$\frac{1}{2b+4}<\frac{1}{\dfrac{a}{1.5}+1},\,\frac{1}{\dfrac{b}{1.5}+1}$$
and
$$\frac{(k_{i-1}+1)p_{i-1}}{1.5}+1<\frac{(b+1.5)p_{i-1}}{1.5}+1<2(b+1.5)+1=2b+4,$$
for $i=1,\ldots,n$. So, we just have checked the condition (\ref{eq15}) for all the generating points $a$, $b$ and $(k_{i-1}+1)p_{i-1}$, $i=1,\ldots,n$.

Now, by the setting (\ref{eq26}), when $k_n=n_0$ with $n\ge 3$, we have constructed the piecewise function $F$ that approximates the sine function on
$[a,b]$ with $0.8\le a$, depending on where $b$ is from $(k_{n-1} + 1)p_{n-1}$ (see Figure \ref{Fig5}):
\begin{figure}[htbp]
\centering
\includegraphics{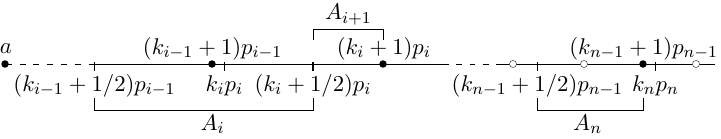}
\caption{In the case of $k_n=n_0$: the black nodes indicate generating points, and the white nodes indicate the possible positions of $b$.}\label{Fig5}
\end{figure}
\begin{itemize}
\item $b\le(k_{n-1}+1/2)p_{n-1}$:
\begin{equation}\label{eq28}
F=\left\{\begin{array}{@{}lcl}
A_0&&\text{on $[a,(k_0+1/2)p_0)$,}\\
A_i&&\text{on $[(k_{i-1}+1/2)p_{i-1},(k_i+1/2)p_i),\quad i=1,\ldots,n-2$,}\\
A_{n-1}&&\text{on $[(k_{n-2}+1/2)p_{n-2},b]$.}
\end{array}
\right.
\end{equation}
\item $(k_{n-1}+1/2)p_{n-1}<b\le(k_{n-1}+1)p_{n-1}$:
\begin{equation}\label{eq29}
F=\left\{\begin{array}{@{}lcl}
A_0&&\text{on $[a,(k_0+1/2)p_0)$,}\\
A_i&&\text{on $[(k_{i-1}+1/2)p_{i-1},(k_i+1/2)p_i),\quad i=1,\ldots,n-1$,}\\
A_n&&\text{on $[(k_{n-1}+1/2)p_{n-1},b]$.}
\end{array}
\right.
\end{equation}
\item $(k_{n-1}+1)p_{n-1}<b$:
\begin{equation}\label{eq30}
F=\left\{\begin{array}{@{}lcl}
A_0&&\text{on $[a,(k_0+1/2)p_0)$,}\\
A_i&&\text{on $[(k_{i-1}+1/2)p_{i-1},(k_i+1/2)p_i),\quad i=1,\ldots,n-1$,}\\
A_n&&\text{on $[(k_{n-1}+1/2)p_{n-1},(k_{n-1}+1)p_{n-1})$,}\\
B_0&&\text{on $[(k_{n-1}+1)p_{n-1},b]$.}
\end{array}
\right.
\end{equation}
\end{itemize}

We just have constructed the piecewise approximate function $F$ on $[a,b]$ in the case of $0.8\le a$. Here we give a summary of our construction process in Algorithm \ref{HQ-Algo2}.
\begin{algthm}\label{HQ-Algo2}
Determining the piecewise function $F$ on $[a,b]$ from the rational numbers $a,b$ with $0.8\le a$ and an integer $r>0$.
\end{algthm}
\begin{algorithmic}[1]
\REQUIRE Two rational numbers $a,b$ ($0.8\le a<b$) and an integer $r>0$;
\ENSURE The list $[[a_i,b_i],f_i]$, $i=1,\ldots,N$, where $\sin\approx F=f_i$ on $[a_i,b_i]$ with the accuracy of $1/10^r$;
\STATE $m:=r+1$;
\WHILE{$\dfrac{1}{(2b+4)10^{r+1}}<\dfrac{1}{10^m}$}
\STATE $m:=m+1$; \COMMENT{The value of $m$ just outside of the loop will be used with \texttt{IntvApprox} for all generating points below};
\ENDWHILE
\STATE $q_0:=\texttt{IntvApprox($b$)[$1$]}$, $n_0:=\texttt{IntvApprox($b$)[$2$]}$, $B_0:=\texttt{IntvApprox($b$)[$3$]}$;
\STATE $i:=0$;
\STATE $p_i:=\texttt{IntvApprox($a$)[$1$]}$, $k_i:=\texttt{IntvApprox($a$)[$2$]}$, $A_i:=\texttt{IntvApprox($a$)[$3$]}$;
\STATE $u:=k_i$;
\IF{$n_0\le u$}
\RETURN $[[a,b],A_0]$;
\ENDIF
\WHILE{$u<n_0$}
\STATE $i:=i+1$;
\STATE $p_i:=\texttt{IntvApprox($(k_{i-1}+1)p_{i-1}$)[$1$]}$;
\STATE $k_i:=\texttt{IntvApprox($(k_{i-1}+1)p_{i-1}$)[$2$]}$;
\STATE $A_i:=\texttt{IntvApprox($(k_{i-1}+1)p_{i-1}$)[$3$]}$;
\STATE $u:=k_i$;
\ENDWHILE
\IF{$i=1$}
 \IF{$b\le(k_0+1/2)p_0$}
 \RETURN $[[a,b],A_0]$;
 \ELSIF{$b\le(k_0+1)p_0$}
 \RETURN $[[a,(k_0+1/2)p_0],A_0]$, $[[(k_0+1/2)p_0,b],A_1]$;
 \ELSE
 \RETURN $[[a,(k_0+1/2)p_0],A_0]$, $[[(k_0+1/2)p_0,(k_0+1)p_0],A_1]$, $[[(k_0+1)p_0,b],B_0]$;
 \ENDIF
\ELSIF{$i=2$}
 \IF{$b\le(k_1+1/2)p_1$}
 \RETURN $[[a,(k_0+1/2)p_0],A_0],[[(k_0+1/2)p_0,b],A_1]$;
 \ELSIF{$b\le(k_1+1)p_1$}
 \RETURN $[[a,(k_0+1/2)p_0],A_0]$, $[[(k_0+1/2)p_0,(k_1+1/2)p_1],A_1]$, $[(k_1+1/2)p_1,b],A_2]$;
 \ELSE
 \RETURN The list derived from (\ref{eq27}) (similar to the ones above);
 \ENDIF
\ELSE
 \IF{$b\le(k_{i-1}+1/2)p_{i-1}$}
 \RETURN The list derived from (\ref{eq28}) (similar to the ones above);
 \ELSIF{$b\le(k_{i-1}+1)p_{i-1}$}
 \RETURN The list derived from (\ref{eq29}) (similar to the ones above);
 \ELSE
 \RETURN The list derived from (\ref{eq30}) (similar to the ones above);
 \ENDIF
\ENDIF
\vskip1ex
\hrule
\vskip2ex
\end{algorithmic}

Now, in the case of $0\le a<0.8$, we can determine $F$ on $[a,b]$ as follows:
\begin{itemize}
\item If $b\le 0.8$: we set $F=F_0$ on $[a,b]$ and we can write $F:\quad[[a,b],F_0]$.
\item If $0.8<b$: we set $F=F_0$ on $[a,0.8]$ and we determine $F$ on $[0.8,b]$ from Algorithm \ref{HQ-Algo2}. So we can write
$$F:\quad [[a,0.8],F_0],\,[[0.8,b_1],f_1],\,[[a_2,b_2],f_2],\,\ldots,\,[[a_n,b],f_n].$$
\end{itemize}

To sum up, given rational numbers $a,b$ with $0\le a<b$ and an integer $r>0$, we can give the piecewise approximate function $F$ to the sine function
with the accuracy of $1/10^r$ in the form of a list $[[a_i,b_i],f_i]$, $i=1,2,\ldots,N$, where $a_1=a$, $b_N=b$, $\sin\approx F=f_i$ on $[a_i,b_i]$. This result
can be described through a procedure, say \texttt{PiecewiseFunctP}, that takes $a,b,r$ as its input and gives the list as its output, then we can write
$$\texttt{PiecewiseFunctP($a$,$b$,$r$)}=[[a_1,b_1],f_1],\,[[a_2,b_2],f_2],\,\ldots,\,[[a_N,b_N],f_N].$$
We can also access each ordered component of such a list with the aid of current computational tools.

Finally, we give a solution to the initial problem: Find a piecewise function $F$ that approximates the sine function
on an interval $[a,b]$ with the accuracy of $1/10^r$ for any rational numbers $a,b$ ($a<b$). The solution $F$ is determined in the following cases:
\begin{enumerate}
\item If $0\le a$: $F$ on $[a,b]$ is given in the form $F:\quad\texttt{PiecewiseFunctP($a$,$b$,$r$)}$.
\item If $a<b\le0$: Suppose that we have
$$\texttt{PiecewiseFunctP($-b$,$-a$,$r$)}=[[a_1,b_1],A_1],\,[[a_2,b_2],A_2],\,\ldots,\,[[a_n,b_n],A_n].$$
So, the sequence
\begin{equation}\label{eq31}
[-b_n,-a_n],\,[-b_{n-1},-a_{n-1}],\,\ldots,\,[-b_1,-a_1]
\end{equation}
is a partition of $[a,b]$ into approximate intervals. Then we set $F=B_j$ on $[-b_j,-a_j]$ with $B_j(x)=-A_j(-x)$, $j=1,\ldots,n$.
Therefore, in this case, these settings are combined with (\ref{eq31}) to become components of the output of a procedure, say \texttt{PiecewiseFunctN}, that takes $a,b,r$
as its input, and we can write
$$F:\quad\texttt{PiecewiseFunctN($a$,$b$,$r$)}=[[-b_n,-a_n],B_n],\,\ldots,\,[[-b_1,-a_1],B_1].$$
\item If $a<0<b$: Thanks to
 the procedures \texttt{PiecewiseFunctP} and \texttt{PiecewiseFunctN}, $F$ on $[a,b]$ can be determined by
$$F:\quad \texttt{PiecewiseFunctN($a$,$0$,$r$)},\,\texttt{PiecewiseFunctP($0$,$b$,$r$)}.$$
\end{enumerate}

Similarly, we can find a solution $G$ to the problem: Find a piecewise function $G$ that approximates the cosine function
on an interval $[a,b]$ with the accuracy of $1/10^r$ for any rational numbers $a,b$ ($a<b$), noting that the cosine function is even.

We will give some numerical and graphical examples of exploitation of our algorithms thanks to the choice of suitable tools. Any system that can do both symbolic and numeric
computation may be a good choice. Here we use Maple to do our calculations. We can refer to \cite{three,four} and Maple help pages for
details about meaning, syntax and usage of Maple commands. 

\begin{figure}[htbp]
\centering
\includegraphics[width=7cm]{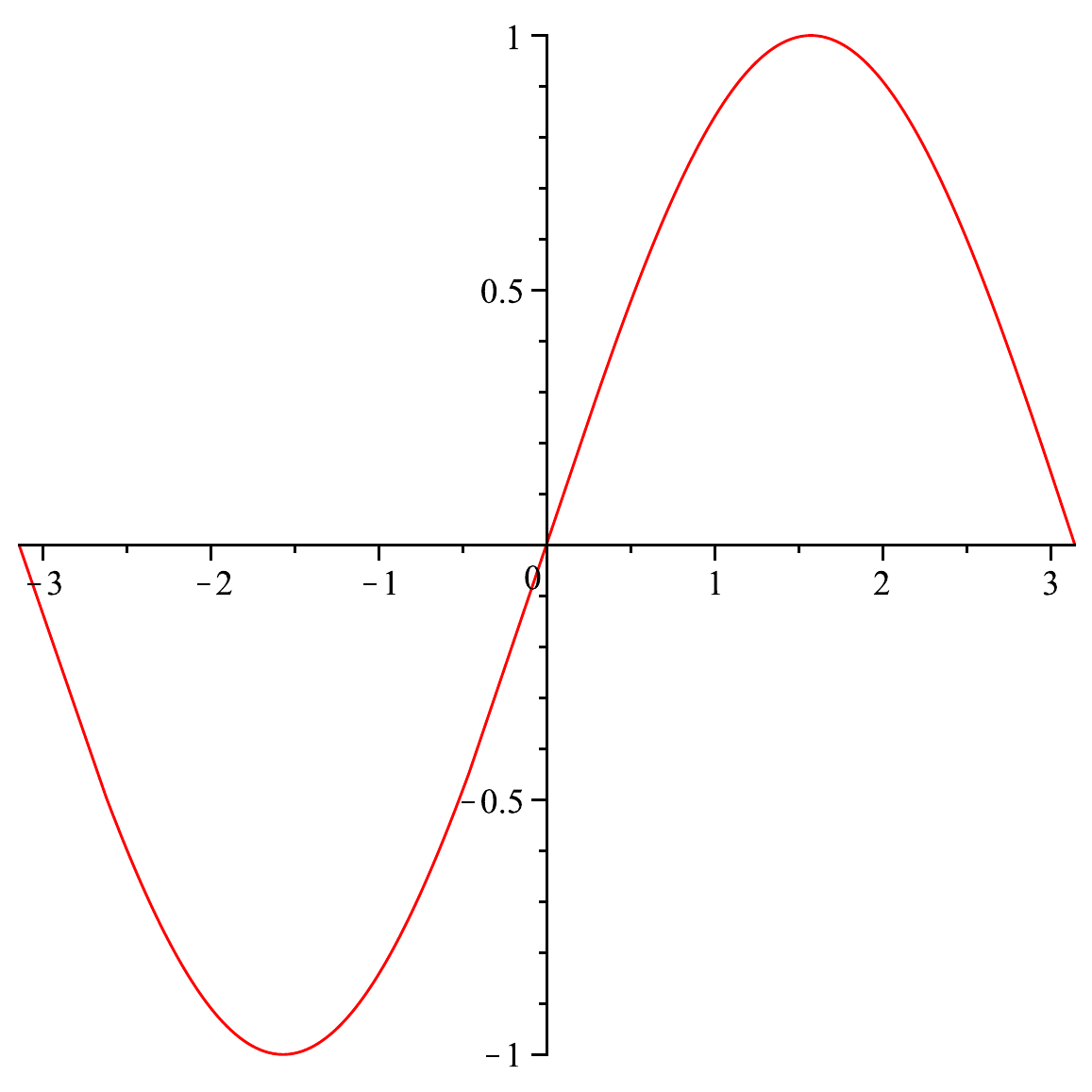}
\caption{The graph of $F\approx \sin$ with $r=12$.}\label{Fig6}
\end{figure}

The first example is for the piecewise function $F$ on $[-3.1416,3.1416]$ such that $\sin\approx F$ on this interval with the accuracy of $1/10^{12}$. We simply write ``with $r=12$'' for such a
precision. Then, the interval can be divided into $8$ approximate intervals as follows
\begin{align*}
&[-3.1416,-3.141592654],\,[-3.141592654,-2.356194490],\,[-2.356194490,-0.8],\,[-0.8, 0],\\
&[0,0.8],\,[0.8,2.356194490],\,[2.356194490,3.141592654],\,[3.141592654,3.1416].
\end{align*}
We can also access all the approximate polynomials of the sine function on these intervals with $r=12$. For instance,
if $\sin\approx F$ on $[2.356194490,3.141592654]$, then $F$ can be the polynomial below:
\begin{align*}
F(x)&=-1.605904384\times10^{-10}x^{13}+6.558626638\times10^{-9}x^{12}-9.857509120\times10^{-8}x^{11}\\
&+5.583443911\times10^{-7}x^{10}-3.414625382\times10^{-7}x^9+0.00001299712087x^8\\
&-0.0002403477863x^7+0.0001044730776x^6+0.008133518439x^5+0.0002885529461x^4\\
&-0.1669715173x^3+0.0002225801191x^2+0.9998995297x+0.00002114256756.
\end{align*}
By taking some points in the intervals, we can check the precision of the approximations to the sine function by the values of $F$ at these points; for instance,
$$\sin(2.5)\approx F(2.5)=0.598472144105;\quad \sin(3.1)\approx F(3.1)=0.0415806624335.$$
Now we give the graph of $F$ on $[-3.1416,3.1416]$ in Figure \ref{Fig6}.

\begin{figure}[htbp]
\centering
\begin{tabular}{@{}cc}
\includegraphics[width=7cm]{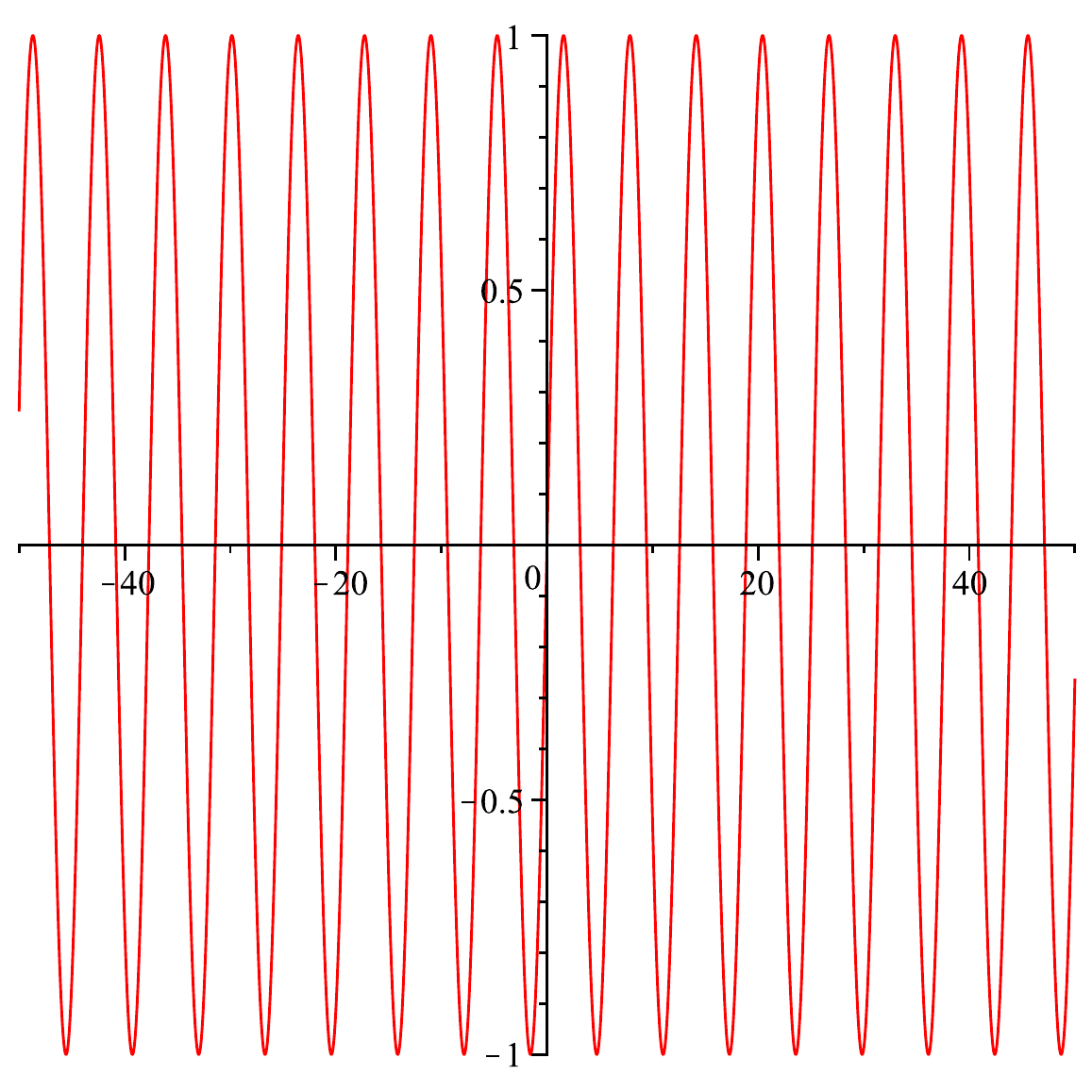}&\includegraphics[width=7cm]{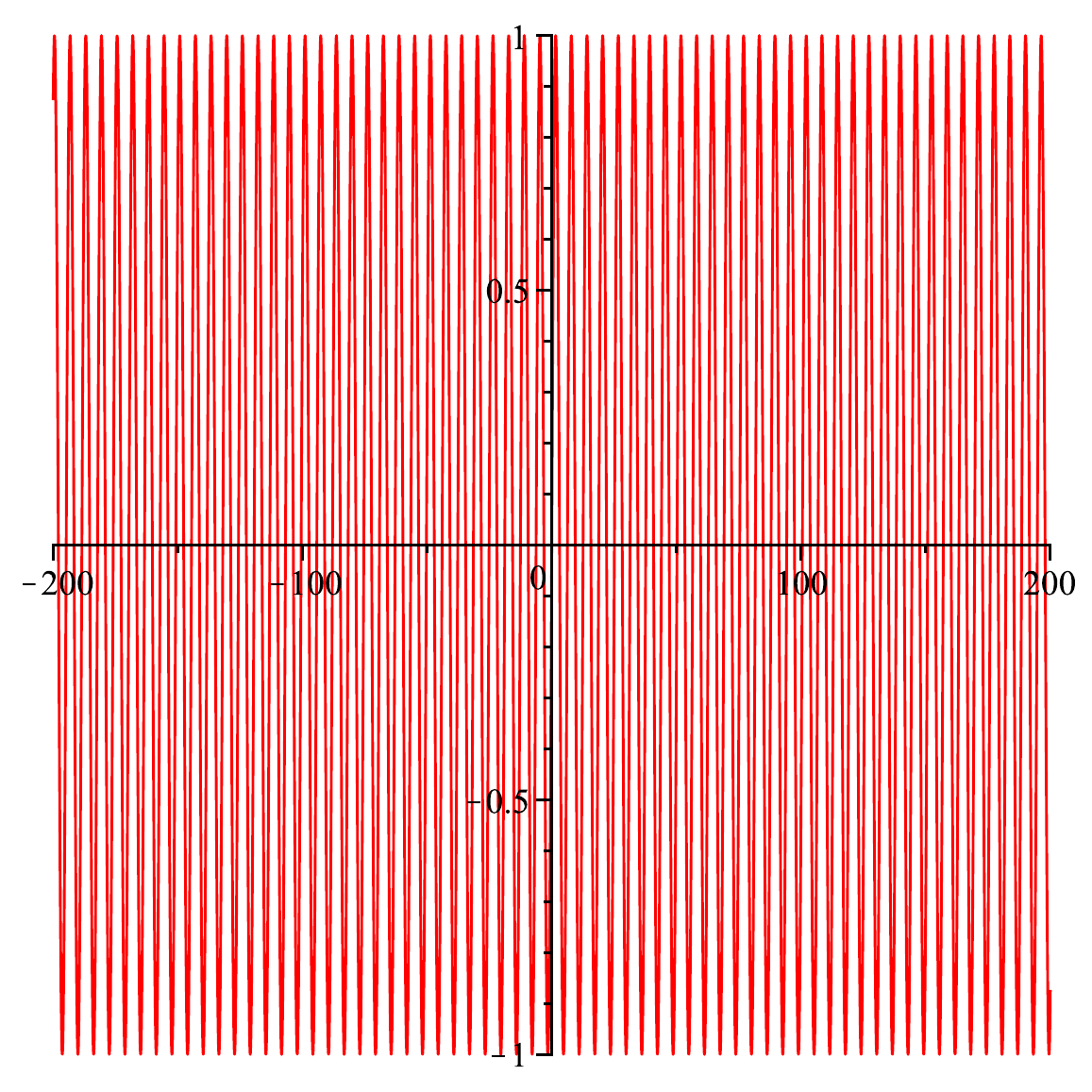}
\end{tabular}
\caption{Left: The graph of $F_1$ with $r=50$. Right: The graph of $F_2$ with $r=200$.}\label{Fig7}
\end{figure}

The second example is for the piecewise functions $F_1$ on $[-50,50]$ with $r=50$ and $F_2$ on $[-200,200]$ with $r=200$.
The intervals $[-50,50]$ and $[-200,200]$ can be divided into $66$ and $258$ approximate intervals, respectively.
First we can check the precision of the approximation below:
$$\sin(49)\approx F_1(49)=-0.95375265275947181836042355858771059528293218973123\quad(\text{with $r=50$}).$$
Then, we give the graphs of $F_1$ and $F_2$ in Figure \ref{Fig7}.

\section{Some other approximation schemes}\label{sect4}

It would be necessary to recall some trigonometric formulas that may be used to make other approximation algorithms. Here they are:
\begin{equation}\label{eq32}
\sin y=3\sin\Big(\frac{y}{3}\Big)-4\sin^3\Big(\frac{y}{3}\Big),\quad\cos y=4\cos^3\Big(\frac{y}{3}\Big)-3\cos\Big(\frac{y}{3}\Big)
\end{equation}
and
\begin{equation}\label{eq33}
\sin y=\sin\Big(\frac{y}{3}\Big)\Big[1+2\cos\Big(\frac{2y}{3}\Big)\Big],\quad \cos y=\cos\Big(\frac{y}{3}\Big)\Big[2\cos\Big(\frac{2y}{3}\Big)-1\Big].
\end{equation}

For $|y|<1$, we may have the following approximations to the sine function:
\begin{equation}\label{eq34}
\sin y\approx3P_n\Big(\frac{y}{3}\Big)-4\Big[P_n\Big(\frac{y}{3}\Big)\Big]^3
\end{equation}
and
\begin{equation}\label{eq35}
\sin y\approx P_n\Big(\frac{y}{3}\Big)\Big[1+2Q_n\Big(\frac{2y}{3}\Big)\Big].
\end{equation}
We try some estimates for their absolute errors. We first consider the approximation (\ref{eq34}). Putting $\sin(y/3)=P_n(y/3)+\varepsilon_s$, we can have the estimates
\begin{equation}\label{eq36}
\Big|P_n\Big(\frac{y}{3}\Big)\Big|\le\Big|\sin\Big(\frac{y}{3}\Big)\Big|+|\varepsilon_s|\le\frac{|y|}{3}+|\varepsilon_s|,\quad\text{where $|\varepsilon_s|\le\dfrac{|y|^{n+1}}{3^{n+1}(n+1)!}$.}
\end{equation}
Then, from (\ref{eq32}), we can write the error $\varepsilon$ of (\ref{eq34}) in the form
$$\varepsilon=3\varepsilon_s-4\varepsilon_s^3-12\varepsilon_sP_n\Big(\frac{y}{3}\Big)\sin\Big(\frac{y}{3}\Big).$$
By the estimates in (\ref{eq36}), we obtain an estimate for $|\varepsilon|$ as
\begin{equation}\label{eq37}
|\varepsilon|\le 4\Big(\frac{|y|^{3n+3}}{[3^{n+1}(n+1)!]^3}+\frac{|y|^{2n+3}}{[3^{n+1}(n+1)!]^2}\Big)+[(4/3)y^2+3]\frac{|y|^{n+1}}{3^{n+1}(n+1)!}.
\end{equation}
The right-hand side of (\ref{eq37}) is denoted by $\varepsilon_1(y,n)$.

Now we consider the approximation (\ref{eq35}). Putting $\cos(2y/3)=Q_n(2y/3)+\varepsilon_c$, we have the estimate
\begin{equation}\label{eq38}
|\varepsilon_c|\le\Big(\frac{2}{3}\Big)^{n+1}\frac{|y|^{n+1}}{(n+1)!}.
\end{equation}
Then, from (\ref{eq33}), we can write the error $\varepsilon$ of (\ref{eq35}) in the form
$$\varepsilon=\varepsilon_s+2P_n\Big(\frac{y}{3}\Big)\varepsilon_c+2\varepsilon_s\cos\Big(\frac{2y}{3}\Big).$$
By the estimates in (\ref{eq36}) and (\ref{eq38}), we obtain an estimate for $|\varepsilon|$ as
\begin{equation}\label{eq39}
|\varepsilon|\le\frac{1}{3^n}\Big(\frac{2}{3}\Big)^{n+2}\frac{|y|^{2n+2}}{[(n+1)!]^2}+\frac{|y|^{n+1}}{3^n(n+1)!}+\Big(\frac{2}{3}\Big)^{n+2}\frac{|y|^{n+2}}{(n+1)!}.
\end{equation}
The right-hand side of (\ref{eq39}) is denoted by $\varepsilon_2(y,n)$.

We can see that the dominant terms in $\varepsilon_1(y,n)$ and $\varepsilon_2(y,n)$ are respectively
$$[(4/3)y^2+3]\frac{|y|^{n+1}}{3^{n+1}(n+1)!}\quad\text{and}\quad \Big(\frac{2}{3}\Big)^{n+2}\frac{|y|^{n+2}}{(n+1)!},$$
and both of them are less than $\dfrac{|y|^{n+1}}{(n+1)!}$ when $|y|<1$. However, if we choose $n$ such that
$$\varepsilon_1(0.8,n)<\frac{1}{10^{r+1}}\quad\text{or}\quad\varepsilon_2(0.8,n)<\frac{1}{10^{r+1}},$$
the polynomial of (\ref{eq34}) or (\ref{eq35}) could have the degree of $3n$ or $2n$, respectively. Denoting $\dfrac{(0.8)^{n+1}}{(n+1)!}$ by $\varepsilon_0(n)$, we
can give the degree of approximate polynomials of (\ref{eq34}), (\ref{eq35}) and of our piecewise algorithm, corresponding to the given precision, as in the following table:
\begin{center}
\begin{tabular}{@{}|c|c|c|c|c|c|@{}}\hline
&\multicolumn{2}{|c|}{$\varepsilon_1(0.8,n)<10^{-(r+1)}$}&\multicolumn{2}{|c|}{$\varepsilon_2(0.8,n)<10^{-(r+1)}$}&$\varepsilon_0(n)<10^{-(r+1)}$\\ \cline{2-6}
\raisebox{1.5ex}[0pt][0pt]{$r$}&$n$&$3n$&$n$&$2n$&$n$\\ \hline
$10$&$9$&$27$&$11$&$22$&$12$\\
$20$&$15$&$45$&$18$&$36$&$20$\\
$50$&$30$&$90$&$35$&$70$&$39$\\
$100$&$53$&$159$&$61$&$122$&$66$\\
$200$&$94$&$282$&$106$&$212$&$115$\\ \hline
\end{tabular}
\end{center}

Similarly, for $|y|<1$, we may have the following approximations to the cosine function:
\begin{equation}\label{eq40}
\cos y\approx4\Big[Q_n\Big(\frac{y}{3}\Big)\Big]^3-3Q_n\Big(\frac{y}{3}\Big)
\end{equation}
and
$$\cos y\approx Q_n\Big(\frac{y}{3}\Big)\Big[2Q_n\Big(\frac{2y}{3}\Big)-1\Big].$$
Putting $\cos(y/3)=Q_n(y/3)+\varepsilon_c$ and taking the estimate
$$\Big|Q_n\Big(\frac{y}{3}\Big)\Big|\le 1+|\varepsilon_c|,\quad\text{where $|\varepsilon_c|\le\dfrac{|y|^{n+1}}{3^{n+1}(n+1)!}$,}$$
we can have the following estimate for the absolute error $|\omega|$ of (\ref{eq40}):
\begin{equation}\label{eq41}
|\omega|\le 4\Big(\frac{|y|^{3n+3}}{[3^{n+1}(n+1)!]^3}+3\frac{|y|^{2n+2}}{[3^{n+1}(n+1)!]^2}\Big)+5\frac{|y|^{n+1}}{3^n(n+1)!}.
\end{equation}
The right-hand side of (\ref{eq41}) is denoted by $\omega_1(y,n)$.

In fact, we can use (\ref{eq34}) and (\ref{eq40}) to make another pointwise approximation algorithm for the sine function. By using Algorithm \ref{HQ-Algo1}, but changing the choice of $n$, we can approximate
the value $\sin y$ in the following cases:
\begin{enumerate}
\item $k_0$ is even: We determine $n$ such that
$$\varepsilon_1(y-k_0p',n)<\frac{1}{10^{r+1}},$$
then we set
$$A=P_n\Big(\frac{y-k_0p'}{3}\Big),\quad \sin y\approx(-1)^{k_0/2}(3A-4A^3).$$
\item $k_0$ is odd: We determine $n$ such that
$$\omega_1(y-k_0p',n)<\frac{1}{10^{r+1}},$$
then we set
$$B=Q_n\Big(\frac{y-k_0p'}{3}\Big),\quad \sin y\approx(-1)^{(k_0-1)/2} (4B^3-3B).$$
\end{enumerate}
The advantage of these settings is that both $A$ and $B$ are evaluated from the polynomials whose degrees are smaller than those of the polynomials in
our pointwise algorithm.

Finally, we also recall that for the cosine function, after applying Algorithm \ref{HQ-Algo1} to $y$, we obtain $p'$, $k_0$ and $n$, and we have
$$\cos y=(-1)^{k_0/2}\cos(y-k_0p)\quad\text{or}\quad \cos y=(-1)^{(k_0+1)/2}\sin(y-k_0p),$$
depending on $k_0$ is even or odd. Then, we can respectively approximate the value $\cos y$ with the accuracy of $1/10^r$ by
$$A^*=A^*(y)=(-1)^{k_0/2}Q_n(y-k_0p')\quad\text{or}\quad B^*=B^*(y)=(-1)^{(k_0+1)/2}P_n(y-k_0p').$$

\section{Conclusion}\label{sect5}

The reader can refer to \cite{six,seven} for the error estimates when using our piecewise approximation algorithm to evaluate definite integrals involving powers of trigonometric functions, as well as
to find a desired approximation on $[a,b]$ in $L^2$-norm to the best approximation of the sine function in $\mathcal{P}_{\ell}$, the vector space of polynomials of degree at most $\ell$.
We expect that some application results of our algorithms would supply some more improved computational tools in applied mathematics.

Moreover, we hope that Algorithm \ref{HQ-Algo1} and Algorithm \ref{HQ-Algo2} would play a significant role in seeking new polynomial approximations to the trigonometric functions.

\section*{Acknowledgements}
The author would like to thank the Maplesoft experts for their great work in developing Maple, a powerful and user-friendly product.

\end{document}